\documentclass[a4paper,english,10pt,oneside]{article}\def\zibreport{1}
\def\arxiv{1}

\usepackage{ifthen}

\usepackage{type1cm}        
\usepackage{makeidx}         
\usepackage{graphicx}        
\usepackage{multicol}        
\usepackage[bottom]{footmisc}

\usepackage{newtxtext}       %
\usepackage[varvw]{newtxmath}       

\makeindex             

\usepackage{todonotes}
\usepackage{hyperref}
\usepackage{cleveref}

\usepackage{algorithm}
\usepackage{algorithmic}

\newcommand{\TheTitle}{A GPU accelerated variant of Schroeppel-Shamir's algorithm for solving the market split problem}
\newcommand{\TheAuthors}{N.--C. Kempke, T. Koch}

\hypersetup{%
    pdftitle={\TheTitle},
    pdfauthor={\TheAuthors}
}

\ifthenelse{\zibreport = 1}{
    \usepackage{./zibtitlepage}
    \let\svhline\hline
    \providecommand{\Description}[1]{}
    \newenvironment{acknowledgement}{
    \section*{Acknowledgements}
    }{}
}{}

\begin{document}

\ifthenelse{\zibreport = 1}{
    \ZTPTitle{\TheTitle}
    \title{\TheTitle}

    \ZTPAuthor{
        \ZTPHasOrcid{Nils--Christian Kempke}{0000-0003-4492-9818},
        \ZTPHasOrcid{Thorsten Koch}{0000-0002-1967-0077}
    }
    \author{
        \ZTPHasOrcid{Nils--Christian Kempke}{0000-0003-4492-9818},\and\
        \ZTPHasOrcid{Thorsten Koch}{0000-0002-1967-0077}
    }

    \ZTPInfo{Preprint}
    \ZTPNumber{25-10}
    \ZTPMonth{February}
    \ZTPYear{2026}

    \date{\normalsize February 20, 2026}
    \ifthenelse{\arxiv = 0}{
        \zibtitlepage
    }{}
    \maketitle
}{
    \title*{{\TheTitle}}
    \titlerunning{GPU accelerated Schroeppel-Shamir for solving the market split problem}
    
    \author{Nils-Christian Kempke\orcidID{0000-0003-4492-9818} and\\ Thorsten   Koch\orcidID{0000-0002-1967-0077}}

    \institute{Nils-Christian Kempke \at Zuse Institute Berlin, Berlin, Germany \email{kempke@zib.de}
    \and Thorsten Koch \at Zuse Institute Berlin, Berlin, Germany \\
    Technische Universität Berlin, Berlin, Germany \email{koch@zib.de}}
}

%
%
\maketitle

\abstract{The market split problem (MSP), introduced by Cornuéjols and Dawande (1998), is a challenging binary optimization problem on which state-of-the-art linear programming-based branch-and-cut solvers perform poorly. We present a novel algorithm for solving the feasibility version of this problem, derived from Schroeppel–Shamir's algorithm for the one-dimensional subset sum problem. Our approach is based on exhaustively enumerating one-dimensional solutions of MSP and utilizing GPUs to evaluate candidate solutions across the entire problem. The resulting hybrid CPU-GPU implementation significantly outperforms a parallel CPU-only variant, efficiently solving instances with up to 10 constraints and 90 variables. We demonstrate the algorithm's performance on benchmark problems, solving instances of size (9, 80) in less than fifteen minutes and (10, 90) in up to one day. Given our results, sorting based algorithms can be considered competitive for solving the MSP on modern hardware.}

\section{Introduction}\label{sec:Introduction}

The \emph{market split problem} (MSP) as in \cite{CournejolsDawande_AClassOfHardSmall01Programs}, is given as the optimization problem
\begin{eqnarray*}
    \min & \sum_{i=1}^m |s_i| &\\
    \text{s.t.} & \sum_{j=1}^n a_{ij} x_j + s_i = d_i, &i=1,\dots,m \\
    & x_j \in \{0,1\}, &j=1,\dots,n \\
    & s_i \in \mathbb{Z}, &i=1,\dots,m.
\end{eqnarray*}
Here, $x_{j}$ are binary decision variables, $m, n \in \mathbb{N}$, and we assume $a_{ij}, d_i \in \mathbb{N}_0$.
The feasibility version of MSP (fMSP) is equivalent to the $n$--dimensional subset sum problem ($n$--SSP): Find a vector $x \in \{0,1\}^n$ such that
\begin{equation}\label{eq:nSSP}
    \sum_{j=1}^n a_{ij} x_j = d_i \quad i = 1,\dots,m.
\end{equation}
%
The $n$--SSP is NP-complete \cite{Karp_ReducibilityAmongCombinatorialProblems}.
As in \cite{CournejolsDawande_AClassOfHardSmall01Programs}, $n$--SSP can be reduced to $1$--SSP by introducing a \emph{surrogate constraint}: Given $D\in \mathbb{N}, D > a_{ij}$, we replace the set of equations in \cref{eq:nSSP} with
\begin{equation}\label{eq:surrogate}
	\sum_{i=1}^m (nD)^{i-1} \sum_{j=1}^n a_{ij}x_j = \sum_{i=1}^m (nD)^{i-1} d_i
\end{equation}
and obtain an equivalent $1$--SSP. Following \cite{CournejolsDawande_AClassOfHardSmall01Programs}, this article examines $n$--SSPs where, given $m\in \mathbb{N}$ we set $n = 10 (m-1)$, $a_{ij}$ is chosen uniformly in the range $[0,99]$ and $b_i = \lfloor \frac{1}{2} \sum_{j=1}^n a_{ij}\rfloor$. These problems are often referred to as $(m,n)$ for a given $m$. In \cite{AardalEtAl_1999_MarketSplitAndBasisReduction} the authors show that this choice of $m$ and $n$ leads to a set of hard $n$--SSPs exhibiting very few expected solutions ranging from $0.19$ expected solutions for $m=4$ to $3.32$ for $m=8$, monotonically growing with increasing $m$.

Several techniques for solving fMSP have been proposed. Branch and cut and branch and bound performs particularly poorly on these instances \cite{AardalEtAl_1999_MarketSplitAndBasisReduction,CournejolsDawande_AClassOfHardSmall01Programs,Wang_MarketSplitBranchAndCut,Wu_MarketSplitDistributed,Vogel_MarketSplitHeuristic}.
The main reason is the vast number of linear basic solutions with value $0$ and little pruning during tree exploration. Dynamic programming and sorting based methods suggested in \cite{CournejolsDawande_AClassOfHardSmall01Programs} rely on using a surrogate constraint \cref{eq:surrogate} to reduce the problem to $1$--SSP.
Lattice-based reduction techniques have proven most successful for solving $n$-SSP. These approaches include basis reduction with polytope shrinkage \cite{CournejolsDawande_AClassOfHardSmall01Programs}, basis reduction with linear programming \cite{AardalEtAl_1999_MarketSplitAndBasisReduction}, and basis reduction with lattice enumeration \cite{Wassermann2002_AttackingTheMarketSplitProblem}.

The core contribution of this work is a novel hybrid CPU-GPU co-design that efficiently decomposes and parallelizes the search space, going beyond a simple GPU port to provide a scalable architecture for exact combinatorial search. In this paper, we present a novel GPU-accelerated sorting-based method for solving $n$--SSP not relying on a surrogate constraint. Our approach solves instances up to $(9, 80)$ and $(10,90)$ in often less than fifteen minutes or one day, respectively. We report the to-date fastest times 
for these instances as published in the literature and show that sorting-based methods are competitive for solving $n$--SSP. We compare our implementation to a multi-threaded CPU only variant and demonstrate significant speed-ups of our hybrid algorithm.
%
\\ \textit{Notation} For a given set $S$, we use $2^S := \{s \mid s\subset A\}$ to denote its power set. The cardinality of $2^S$ is given by $2^{|S|}$. We use $\mathbb{N}$ to denote the natural numbers and define $\mathbb{N}_0 := \mathbb{N}\cup \{0\}$.

\section{Enumerating solutions of $1$-SSP}\label{sec:1dssp}

Our algorithm is based on the exhaustive enumeration of all solutions of $1$--SSP. Let $S := \{1,\dots,n\}$ be an index set for $n \in \mathbb{N}$ and for each $i\in S$ let $a_i \in \mathbb{N}_0$ be an integer weight. For any subset $s \subset S$ we define its subset sum as $a(s) := \sum_{i \in s} a_i$. A classic way to find a solution to $1$--SSP is Horowitz-Sahni's two-list algorithm \cite{HorowitzSahni_ComputingPartitions}, also used in \cite{CournejolsDawande_AClassOfHardSmall01Programs} in combination with a surrogate constraint. The two-list algorithm splits the coefficients of $1$--SSP into two subsets, generates all subset sums of each subset sorted by value, and then traverses the two sorted lists in ascending and descending order until a pair is found whose combined value satisfies the $1$--SSP.
\begin{algorithm}
    \caption{Schroeppel-Shamir's algorithm for all solutions of $1$--SSP} \label{alg:SS}
    \begin{algorithmic}[1]
        \REQUIRE $S = \{1,\dots,n\}$, weights $(a_1,\dots,a_n)$, $a_i \in \mathbb{N}_0$, $d \in \mathbb{N}_0$, weight function $a: 2^S \to \mathbb{N}_0$
        \STATE Initialize set of solutions $\mathcal{R} \gets \emptyset$
        \STATE Partition $S$ into $A, B, C, D$ of size $\approx\frac{n}{4}$
        \STATE Let $(s^A_1, \dots, s^A_{2^{|A|}})$ be an ordering of $2^A$ with $a(s^A_1) \leq \dots \leq a(s^A_{2^{|A|}})$
        \STATE Let $(s^C_1, \dots, s^C_{2^{|C|}})$ be an ordering of $2^C$ with $a(s^C_1) \geq \dots \geq a(s^C_{2^{|C|}})$
        \STATE Initialize min-heap $H_1 \gets \{(s^A_1, s^B) \mid s^B \in 2^B\}$ with key $a(s^A_1) + a(s^B)$
        \STATE Initialize max-heap $H_2 \gets \{(s^C_1, s^D) \mid s^D \in 2^D \}$ with key $a(s^C_1) + a(s^D)$
        \WHILE{$H_1$ and $H_2$ not empty}
        \STATE $(s^A_i, s^B) \gets$ top($H_1$), $(s^C_j, s^D) \gets$ top($H_2$)
        \STATE $\alpha := a(s^A_i) + a(s^B)$,  $\beta := a(s^C_j) + a(s^D)$
        \IF{$\alpha + \beta < d$ \AND $i + 1 \leq 2^{|A|}$}
            \STATE $\text{pop}(H_1)$ and insert $(s^A_{i+1}, s^B)$ into $H_1$
        \ELSIF{$\alpha + \beta > d$ \AND $j + 1 \leq 2^{|C|}$}
            \STATE $\text{pop}(H_2)$ and insert $(s^C_{j+1}, s^D)$ into $H_2$
        \ELSE
            \STATE $E_A := \{s \in 2^A \mid a(s) = a(s^A_i)\}$, $E_C := \{s \in 2^C \mid a(s) = a(s^C_j)\}$ \label{alg:SS:line:extract}
            \STATE $\mathcal{R} \gets \mathcal{R} \cup \{s^A \cup s^B \cup s^C \cup s^D \mid s^A \in E_A, s^C \in E_C\}$
            \STATE $\text{pop}(H_1)$
            \STATE \textbf{if} $i+|E_A| \leq 2^{|A|}$ \textbf{then} insert $(s^A_{i+|E_A|}, s^B)$ into $H_1$ \textbf{end if}
            \STATE $\text{pop}(H_2)$
            \STATE \textbf{if} $j+|E_C| \leq 2^{|C|}$ \textbf{then} insert $(s^C_{j+|E_C|}, s^D)$ into $H_2$ \textbf{end if}
        \ENDIF
        \ENDWHILE
    \RETURN $\mathcal{R}$
    \end{algorithmic}
\end{algorithm}
For large $1$--SSP instances, the space complexity $\mathcal{O}(2^\frac{n}{2})$ of Horowitz-Sahni's algorithm quickly becomes prohibitive.

An improvement in space complexity is provided by the algorithm of Schroeppel-Shamir \cite{SchroeppelShamir_AlgorithmForCertainNPCompleteProblems} as shown in \Cref{alg:SS}. Instead of generating the two power sets used by Horowitz-Sahni, it uses heaps to dynamically generate each power set, reducing the space complexity to $\mathcal{O}(2^\frac{n}{4})$. In \Cref{alg:SS}, we modified the original algorithm to collect all solutions of $1$--SSP. The subsets $E_A$ and $E_C$ in line~\ref{alg:SS:line:extract} can be determined quickly by linearly iterating the sorted $2^A$ and $2^C$ staring at the $i$-th and $j$-th elements, respectively.

\section{GPU accelerated Schroeppel-Shamir for the $n$-SSP}\label{sec:GPU_nSSP}

Given \Cref{alg:SS} for retrieving all solutions of the $1$--SSP, we solve $n$-SSP in the following way: For a given $n$--SSP, we use Schroeppel–Shamir's algorithm to find all solutions of the $1$--SSP obtained by only considering the first row of \cref{eq:nSSP}. We note that, in order to be correct, the $1$--SSP algorithm used for $n$--SSP needs to take generate all solutions to $1$-SSP explicitly including zero coefficients as well (as done in \Cref{alg:SS}). Whenever a set of solutions to $1$--SSP is found, we verify it against the rest of the problem. The procedure is shown in \Cref{alg:SS-nSSP}.
\begin{algorithm}[ht]
    \caption{Schroeppel–Shamir for $n$--SSP}\label{alg:SS-nSSP}
    \begin{algorithmic}[1]
        \REQUIRE $S := \{1,\dots,n\}$, $A \in \mathbb{N}_0^{m\times n} = (a_{ij})$, $d\in\mathbb{N}_0^m$
        \STATE Set up heaps $H_1$, $H_2$ as in \Cref{alg:SS} w.r.t. the weights $(a_{1i},\dots,a_{1n})$ and $d_1$.
        \WHILE{$H_1$, $H_2$ not empty}
            \STATE Iterate $H_1$, $H_2$ as before, using $\alpha$, $\beta$, $s^B$, and $s^D$
            \IF{$\alpha + \beta = d_1$}
                \STATE Extract the sets of equal weight $E_A$ and $E_C$ as before
                \label{alg:SS-nSSP:line:collect}
                \FORALL{$s^A\in E_A$, $s^C\in E_C$} \label{alph:SS-nSSP:line:first_val_loop}
                    \STATE Let $x \in \{0,1\}^n$ be the characteristic vector of $s^A \cup s^B \cup s^C \cup s^D$
                    \STATE \textbf{if} $Ax = d$ \textbf{then} \textbf{return} $x$ \textbf{end if}
                \ENDFOR \label{alph:SS-nSSP:line:last_val_loop}
            \ENDIF
        \ENDWHILE
    \end{algorithmic}
\end{algorithm}
To make this approach feasible, we rely on GPU acceleration for the validation loop line~\ref{alph:SS-nSSP:line:first_val_loop} to line~\ref{alph:SS-nSSP:line:last_val_loop} in \Cref{alg:SS-nSSP}. After collecting all solutions in line~\ref{alg:SS-nSSP:line:collect}, we offload the solution validation to the GPU while the CPU continues searching for $1$--SSP solutions. This creates a CPU-GPU hybrid pipeline interleaving collection and validation operations. The GPU validation algorithm is shown in \Cref{alg:gpu}.
\begin{algorithm}[ht]
    \caption{GPU-accelerated solution validation} \label{alg:gpu}
    \begin{algorithmic}[1]
        \REQUIRE $E_A$, $E_C$, $s^B$, and $s^D$ as in \Cref{alg:SS-nSSP} line~\ref{alph:SS-nSSP:line:first_val_loop}; $A$, $b$ as in \Cref{alg:SS-nSSP}
        \STATE Compute $Q := \{Ax \mid x \text{ is characteristic vector of } s^A \cup s^B\,\ s^A \in E_A\}$ \label{alg:gpu:comp1}
        \STATE Compute $R := \{b - Ax \mid x \text{ is characteristic vector of } s^C \cup s^D,\ s^C\in E_C\}$ \label{alg:gpu:comp2}
        \STATE $Q_{\text{enc}} \gets \texttt{encode\_kernel}(Q)$
        \STATE $R_{\text{enc}} \gets \texttt{encode\_kernel}(R)$
        \STATE \texttt{sort\_kernel}$(Q_{\text{enc}})$
        \STATE \texttt{parallel\_binary\_search\_kernel}$(Q_{\text{enc}}, R_{\text{enc}})$
    \end{algorithmic}
\end{algorithm}
Instead of naively checking each pair in the validation loop, we hash the partial subset sum vectors in the \texttt{encode\_kernel}, sort one hash set, and perform a parallel binary search for elements of the unsorted set. For lines~\ref{alg:gpu:comp1} and \ref{alg:gpu:comp2}, we pre-compute and buffer the partial vectors $Ax$ for all characteristic vectors. Sorting and parallel binary search are implemented using CUDA thrust~\footnote{\url{https://nvidia.github.io/cccl/thrust/}}. Encoding uses a custom hash kernel, simplified shown in \Cref{alg:gpu:encode}.
\begin{algorithm}[ht]
    \caption{\texttt{encode\_kernel}: parallel hash encoding} \label{alg:gpu:encode}
    \begin{algorithmic}[1]
        \REQUIRE $d_i \in \mathbb{N}^m$, $i = 0, \dots, N-1$; array \texttt{hash} of size $N$
        \IF{\texttt{threadId} $< N$}
            \STATE $h \gets 0$
            \STATE \textbf{for} $j =0,\dots,m - 1$ \textbf{do} $h \gets \texttt{hash\_two}(h, (d_{\texttt{threadId}})_j)$ \textbf{end for}
            \STATE $\texttt{hash}[\texttt{threadId}] \gets h$
        \ENDIF
    \end{algorithmic}
\end{algorithm}
For large $m$, the number of $1$--SSP solutions passed to \Cref{alg:gpu} grows rapidly, potentially exceeding GPU memory. We then validate $Q$ and $R$ quadratically in chunks by partitioning both arrays. This creates a bottleneck that could be addressed using multiple GPUs, though our current implementation uses a single GPU.

\section{Computational Results}\label{sec:results}

Our implementation is available on GitHub\footnote{\url{https://github.com/NCKempke/MarketShareGpu}}. Experiments were conducted on a NVIDIA GH200 Grace-Hopper super-chip\footnote{\url{https://www.nvidia.com/en-us/data-center/grace-hopper-superchip/}}, with an ARM Neoverse-V2 CPU (72 cores), 480 GB of memory, and one NVIDIA H200 GPU with 96 GB of device memory.
We used the fMSP instances provided in QOBLIB\footnote{\url{https://git.zib.de/qopt/qoblib-quantum-optimization-benchmarking-library/}}~\cite{Koch_QOBLIB}. QOBLIB contains for each $m\in\{3,\dots,15\}$, $K\in\{50,100,200\}$, four feasible fMSP instances with coefficients uniformly drawn in $[0, K)$. We reflect this using the extended notation $(m, n, K)$. We ran each instance with our CPU-GPU hybrid algorithm and, formulated as a linear integer program (as in the original MSP) using Gurobi 11~\cite{Gurobi11}. Additionally, we ran instances for $m\in\{3,\dots,8\}$ with a CPU only parallel implementation of our algorithm relying on OpenMP\footnote{\url{https://www.openmp.org/}} using all 72 available cores (also available on GitHub). We used a time limit of three days.
\begin{table}
\caption{Solutions times in seconds for fMSPs with Schroeppel-Shamir and Gurobi}\label{tab:results}
\begin{tabular}{lrrrrrrr}
\hline\noalign{\smallskip}
\textbf{Class} & Instance 1 & Instance 2 & Instance 3 & Instance 4 & Avg. & Avg. CPU & Avg. Gurobi \\
\noalign{\smallskip}\svhline\noalign{\smallskip}
        $(7, 60, 50)$  & 0.37 & 0.37 & 0.39 & 0.35 & 0.37 & 2.79 & 243.32 \\
        $(7, 60, 100)$  & 1.66 & 0.88 & 1.38 & 0.86 & 1.20 & 19.94 & 1\,086.08 \\
        $(7, 60, 200)$  & 2.07 & 2.45 & 2.35 & 1.19 & 2.02 & 24.15 & 3\,158.37 \\
        $(8, 70, 50)$  & 1.37 & 1.12 & 1.00 & 1.50 & 1.25 & 21.97 & MEM \\
        $(8, 70, 100)$  & 5.80 & 8.07 & 9.19 & 8.28 & 7.84 & 386.21 & MEM \\
        $(8, 70, 200)$  & 15.60 & 20.13 & 8.45 & 27.03 & 17.80 & 879.25 & MEM \\
        $(9, 80, 50)$  & 23.25 & 10.78 & 14.68 & 15.40 & 16.03 & 574.54 & MEM \\
        $(9, 80, 100)$  & 472.88 & 101.52 & 300.37 & 69.51 & 236.07 & 15\,212.71 & MEM \\
        $(9, 80, 200)$  & 548.83 & 505.44 & 486.97 & 541.80 & 520.76 & 34\,391.54 & MEM \\
        $(10, 90, 50)$  & 153.50 & 288.41 & 74.13 & 155.29 & 167.83 & 16\,127.69 & MEM \\
        $(10, 90, 50)^*$ & 2\,957.91 & 2\,234.68 & 3\,866.84 & 2\,144.57 & 2\,803.80 & 203\,164.16 & MEM \\
        $(10, 90, 100)$  & 152\,323.86 & 209\,021.22 & 176\,957.14 & 31\,544.88 & 148\,663.34 & - & MEM \\
        $(10, 90, 100)^*$ & 104\,230.63 & 30\,141.04 & 56\,104.74 & 76\,068.47 & 66\,636.22 &  - & MEM \\
        $(10, 90, 200)$  & - &        - &       -  &        - & - & - & MEM \\
        $(10, 90, 200)^*$ & - &       -  &       -  &       -  & - & - & MEM \\
        $(11,100, 50)$   & 110\,032.34 & 5\,313.42 & 42\,940.19 & 52\,118.21 & 52\,601.04 & - & MEM \\
        $(11,100, 50)^*$  & 34\,151.77 & 34\,669.73 & 54\,455.52 & 44\,320.57 & 41\,399.39 & -& MEM \\
\noalign{\smallskip}\hline\noalign{\smallskip}
\end{tabular}
\Description{Table 1 is titled ``Solution times in seconds for fMSPs with Schroeppel-Shamir and Gurobi''. It lists runtime results for different instance classes of fMSPs. Each row corresponds to a parameter triple like (m,n,K), indicating a specific problem class. The triplets range from (7,60,{50,100,200}) to (10,90,{50,100,200}). The last line is (11,100,50). There are four instances per class, and their individual solution times (in seconds) are listed, followed by the average time for Schroeppel-Shamir and, if available, Gurobi. For smaller instances like (7,60,50), average times are low: around 0.37 seconds for Schroeppel-Shamir and over 240 seconds for Gurobi. As problem size increases, Schroeppel-Shamir's average rises to over 148,000 seconds for (10,90,100).  Gurobi times are only available for (7,60,K) instances. The final entries (11,100,50) report average times of about 50,000s. For (10,90,K) and (11,100,50) the table contains additional rows marked with an asterisk.}
\end{table}
\begin{table}
\caption{Literature results in seconds for fMSPs $(m,n,100)$}\label{tab:literature}
\begin{tabular}{lrrrrrrrrr}
\hline\noalign{\smallskip}
\textbf{Prob. size} & B\&B \cite{CournejolsDawande_AClassOfHardSmall01Programs} & B\&C \cite{CournejolsDawande_AClassOfHardSmall01Programs} & DP & Basis \cite{CournejolsDawande_AClassOfHardSmall01Programs} & Group \cite{CournejolsDawande_AClassOfHardSmall01Programs} & Sort (ours) & LLL \cite{AardalEtAl_1999_MarketSplitAndBasisReduction} & Enum \cite{Wassermann2002_AttackingTheMarketSplitProblem} & Gurobi \\
\noalign{\smallskip}\svhline\noalign{\smallskip}
(3,20)  & 10.3       & 178.67    &  1.11 & 239.40 & 0.12        & 0.17        & -         & 0.01       & 0.13      \\
(4,30)  & 2\,271.65  & -         & 19.67 & -      & 17.96       & 0.20        & -         & 0.08       & 0.78      \\
(5,40)  & -          & -         & -     & -      & 1\,575.58   & 0.22        & 62.4      & 1.01       & 0.81      \\
(6,50)  & -          & -         & -     & -      & 22\,077.32  & 0.29        & 2\,190.0  & 2.16       & 79.09     \\
(7,60)  & -          & -         & -     & -      & -           & 1.20        & -         & 28.2       & 1\,086.08 \\
(8,70)  & -          & -         & -     & -      & -           & 7.84        & -         & 678        & -         \\
(9,80)  & -          & -         & -     & -      & -           & 236.07      & -         & 9\,733     & -         \\
(10,90) & -          & -         & -     & -      & -           & 66\,636.22  & -         & (655\,089) & -         \\
\noalign{\smallskip}\hline\noalign{\smallskip}
\end{tabular}
\Description{Table 2 is titled ``Literature results in seconds for fMSPs (m,n,100)''. It compares the runtime performance of several algorithms across increasing problem sizes. Each row corresponds to a problem size (m,n) with fixed K=100, from (3,20) up to (10,90). Columns report solution times for: B\&B, B\&C, DP, Basis, Group, Sort, LLL, Enum, and Gurobi. Each column has data for a different amount of instances. For small instances (3,20) all methods are fast, with Bas. Enum at 0.01s and Gurobi at 0.13s. Gurobi tops for instances (7,50) with 1086.08 seconds. DP only has data available for (4,30) with 19.67s. Group tops at (6,50) with 22077.32s. Sort tops at (9,80) with 236.07s and (10,90) with 66636.22 seconds. LLL contains data only for (5,40) and (6,50) with 62.4s and 2190.0s, respecitvely. Enum tops at (9,80) with 9733 and (10,90) with 655089s. The last result of the Enum column is enclosed by brackets.
}
\end{table}
We present our results in \cref{tab:results}. The table shows solution times of our CPU-GPU hybrid algorithm per instance, as well as the average runtimes of our CPU-GPU variant, our CPU-only variant, and Gurobi. The CPU-GPU algorithm solves all instances up to $m=10$. For rows where the \textbf{Class} is annotated by ``*'', we first applied a surrogate constraint dimension reduction, reducing the first 2 constraints into one. The reduction shows speed-ups for $(10,90,100)$ decreasing runtime to about 66\,000 seconds. Generally, instances with a smaller coefficient range solve faster, and solution time grows rapidly with increasing $m$. We could not solve the instances $(10, 90, 200)$ or any instance larger than $(11, 100, 50)$. Compared to the parallel CPU implementation, the hybrid implementation is much more competitive and suffers less from an increase in complexity. It outperforms the CPU-only variant by roughly a factor of 10. We could not solve instances larger than (8,70) efficiently on CPU. Gurobi ran out of memory for all instances larger than $(7,60)$\\
A performance analysis of our hybrid algorithm reveals that its scalability is primarily constrained by GPU memory capacity for storing intermediate results and hashes and by CPU single-thread performance for managing the one dimensional search. Further gains would therefore depend more on increased GPU memory bandwidth and capacity than on higher clock speeds. This confirms that the method's advantage stems from exploiting the GPU's architectural strengths, massive parallelism and high memory throughput, rather than raw clock speed alone. On the other hand, the CPU-only variant would strongly benefit from an increased memory bandwidth, as solution validation time dominates the runtime.\\
\Cref{tab:literature} extends the comparison from \cite{Wassermann2002_AttackingTheMarketSplitProblem} to include our approach (note, the basis enumeration result in brackets for $(10,90)$ was obtained with a single instance). The table combines results from different papers run on different (and sometimes quite outdated) hardware. Our algorithm updates the \emph{Sort} column and currently provides the fastest reported results for fMSPs. However, it would be worth re-running all experiments with up-to-date implementations and on state-of-the art hardware to provide a more complete picture, which is out of the scope of this paper.

\section{Conclusion}\label{sec:conclusion}
In this paper, we presented a novel hybrid CPU-GPU approach for solving the $n$--SSP by accelerating a variant of the Schroeppel–Shamir algorithm. Our implementation efficiently solves feasibility market split instances of size up to $(10,90)$ in under one day on average, representing, to our knowledge, the best published results for this problem class. This work demonstrates the potential of heterogeneous computing for exact combinatorial search, where algorithms effectively leveraging both CPU and GPU remain rare. For future work, we plan to explore GPU-accelerated lattice enumeration methods, as described in \cite{Wassermann2002_AttackingTheMarketSplitProblem}, which we believe offer a promising direction for further performance improvements on MSP.

\begin{acknowledgement}
The work for this article has been conducted in the Research Campus MODAL funded by the German Federal Ministry of Education and Research (BMBF) (fund numbers 05M14ZAM, 05M20ZBM, 05M2025).
\end{acknowledgement}


\end{document}